\documentstyle[12pt]{amsart}
\renewcommand{\thetheorem}{\thesection.\arabic{theorem}}

\theoremstyle{remark}

\newcommand{\Cl}{\operatorname{Cl}}
\newcommand{\im}{\operatorname{im}}

\newcommand{\Pic}{\operatorname{Pic}}

\newcommand{\rank}{\operatorname{rank}}
\newcommand{\Hom}{\operatorname{Hom}}

\newcommand{\B}{\operatorname{B}}

\newcommand{\columnspace}{\operatorname{columnspace}}
\newcommand{\diag}{\operatorname{diag}}
\newcommand{\TNemb}{\operatorname{T_Nemb}}
\begin{document}
\title{COMPUTING THE COHOMOLOGICAL BRAUER GROUP OF A TORIC VARIETY}
\author{T. J. Ford}
\address{Department of Mathematics\\Florida Atlantic University\\Boca
Raton, FL 33431}
\email{Ford@@acc.fau.edu}
\thanks{Supported in part by the NSF under grant DMS-9025092. \\ \today}
\maketitle
\medskip

The purpose of this article is to show how one might compute the \'etale
cohomology groups $H^p_{\acute{e}t}(X,G_m)$
in degrees $p=0$, $1$ and $2$
of a toric variety $X$ with coefficients in the sheaf of
units. The method is to reduce the computation down to the problem of
diagonalizing a matrix with integral coefficients.
The procedure outlined in this article has been fully
implemented by the author
as a program written in the ``C'' programming language.

The groups that we want to compute are finitely generated abelian groups.
Our method for computing them is to reduce the problem to a matrix theory
computation involving matrices with integer coefficients.
Let us begin with a brief description of the  algorithms
for matrices with integer coefficients which are to be used.

Let $S$ be an $m \times n$ matrix with integer coefficients. The basic
procedure that we perform on $S$ is ``diagonalization'', or equivalently,
``put $S$ into Smith normal form.'' This means that we find invertible
matrices $X$ and $Y$ so that $XSY = \diag\{d_1, d_2, \ \dots,\ d_s, 0,
\dots,\ 0\}$ where $d_1|d_2|\ \cdots\ |d_s$ and $d_s \not= 0$. Thus $S$ has
column rank $s$.
An algorithm for computing $X$, $Y$ and $d_1, \dots,\ d_s$
can be found in \cite{CR:RTf}.

Once we have computed the matrix $Y$ we can find simultaneous bases for the
column space of $S$ and ${\Bbb Z}^m$. By this we mean a
set of vectors
$\{x_1,
\dots,\ x_s\}$ that extend to a basis for ${\Bbb Z}^m$ and such that a
basis for the column space of $S$ is $\{d_1x_1, d_2x_2, \dots,\ d_sx_s\}$.
We simply take $\{x_1, \dots,\ x_s\}$ to be the columns of the matrix
$X^{-1}\diag\{1,1,\dots,\ 1,0, \dots,\ 0\}$ (where there are $s$ ones)
$= S Y \diag\{d_1^{-1}, d_2^{-1}, \dots,\ d_s^{-1}, 0, \dots,\ 0\}$.
The vectors  $\{x_1, \dots,\ x_s\}$ span the smallest direct summand of
${\Bbb Z}^m$  that contains the column space of $S$. A basis for the kernel
of $S$ exists in the columns numbered $s+1, \dots,\ n$ of the matrix $Y$.

Say the columns of $B$ contain a basis for a direct summand of rank $s$
of ${\Bbb Z}^m$ and $X$ is an invertible matrix such that $X B =
\diag\{1, \dots, \ 1,0, \dots,\ 0\}$. Let A be an $m \times l$ matrix whose
columns are in the column space of $B$. The matrix $X A$ represents the
function which maps the columns of $A$ into the column space of $B$.

Suppose we have an $m \times l$ matrix $A$ and an $n \times m$ matrix
$B$ such that $B A = 0$ and we want to
compute the homology group $\ker{B} / \im{A}$.
First find a basis for the kernel of $B$ as above. Say this basis makes
up the columns of the matrix $K$. Find invertible matrices $X$ and
$Y$ so that $X K Y = \diag\{1, \dots,\ 1\}$. We want
to write the columns of $A$ as linear combinations of the columns of
$K Y$. The matrix for the embedding $\columnspace(A) \to \columnspace(KY)$ is
the $\rank(K) \times l$ matrix $X A$. The
invariant factors of $\ker{B} / \im{A}$ are obtained by diagonalizing this
matrix.

We now discuss how to use the above algorithms to compute the
cohomological Brauer group of a toric variety. Let $N = {\Bbb Z}^r$, $M =
\Hom_{\Bbb Z}(N,{\Bbb Z})$.
Let $\Delta$ denote a fan
on $N_{\Bbb R} =
{\Bbb R}^r$ and $X = \TNemb(\Delta)$ the corresponding toric variety
over the algebraically closed field $k$ of characteristic $0$.
Let $K$ denote the function field of $X$.
The cohomological Brauer group of $X$ is the second \'etale
cohomology group with coefficients in the sheaf of units,
$H^2_{\acute{e}t}(X,G_m)$.
According to \cite[Theorem 1]{DFM:Cbg}, if $\tilde X$ is an equivariant
desingularization of $X$, then there is a split-exact sequence with natural
maps

\begin{equation}
\label{eq1}
0 \rightarrow  H^2_{\acute{e}t}(K/X,G_m) \rightarrow
H^2_{\acute{e}t}(X,G_m) \rightarrow
H^2_{\acute{e}t}(\tilde X,G_m) \rightarrow  0
\end{equation}

\noindent
Sequence (\ref{eq1})
reduces the calculation of $H^2_{\acute{e}t}(X,G_m)$ down to the computation
of the 2 smaller groups in the sequence.
The relative cohomological Brauer group $H^2_{\acute{e}t}(K/X,G_m)$
consists of the 2-cocycles which are generically split.
The group $H^2_{\acute{e}t}(\tilde X,G_m)$
 is naturally isomorphic to the
image of the Brauer group of $X$ in the Brauer group of the function field
of $X$, $\B(K)$.
In order to compute these groups, our method allows us to also compute the
zeroeth  and first degree \'etale cohomology groups without much extra
work. The divisor class group is the Picard group of a suitable open subset
of X and we might as well compute it as well.
The groups which we plan to compute are now enumerated.

\begin{enumerate}
\item The image of the Brauer group of $X$ in $\B(K)$. This group
is naturally isomorphic to the Brauer group $\B(\tilde X) =
H^2_{\acute{e}t}(\tilde X,G_m)$ of an
equivariant
desingularization $\tilde X$ of $X$.
\item The divisor class group $\Cl(X)$.
\item The group of units $H^0_{\acute{e}t}(X,G_m)$.
\item The Picard group  $\Pic(X) = H^1_{\acute{e}t}(X,G_m)$.
\item The relative cohomological Brauer group, $H^2(K/X,G_m)$.
\end{enumerate}

Let $\{\rho_1, \dots,\ \rho_n\} = \Delta(1)$ be the cones in $\Delta$ of
dimension $1$. For each $\rho_i \in \Delta(1)$ choose a
primitive generator $\eta_i \in {\Bbb Z}^r$ such that
${\rho}_i = {\Bbb R}_{{\ge}0} \cdot \eta_i$.

Let $\tilde{X}$ be an equivariant desingularization of $X$. Then
$\tilde{X}$ is the toric variety associated to a fan
$\Delta '$ obtained by
subdividing $\Delta$.
Any maximal cone $\tau \in \Delta '$ is
contained in some maximal cone $\sigma \in \Delta$.
For such a pair $(\tau,\sigma)$,  the group generated by
$\tau \cap N$ is equal to the smallest direct summand of $N$ that contains
$\{\rho \in \Delta(1) | \rho \subseteq \sigma \}$.
Let $N'$ be the group generated by
$\cup_{\tau \in \Delta '} \sigma \cap N$.
The invariant factors of $N/N'$ can be computed using the methods described
above. For each maximal cone $\sigma$ of $\Delta$, find a basis
$L(\sigma)$ for the
smallest direct summand of $N$ containing $\{ \eta_i | \rho_i \subseteq
\sigma \}$.
Set up a matrix whose columns are the elements of $\cup_{\sigma
\in \Delta} L(\sigma)$ and compute the invariant factors of this matrix
using the above mentioned methods. Say the invariants are $a_1,
\dots,\ a_r$.
In order to compute group 1 above, we use \cite[Theorem 1.1]{DF:Bgt} which
says
$H^2_{\acute{e}t}(\tilde X,G_m)$ is isomorphic to
$\oplus_{i=1}^{r-1} \Hom({\Bbb Z}/a_i,{\Bbb Q}/{\Bbb Z})^{r-i}$.
Note that $\Hom({\Bbb Z}/a_i,{\Bbb Q}/{\Bbb Z})= 0, {\Bbb Z}/a_i,
{\Bbb Q}/{\Bbb Z}$ according to  whether $|a_i| = 1, |a_i| > 1, a_i = 0$.

The divisor class group of $X$ is isomorphic to that of the toric variety
associated to the fan $\{0, \rho_1, \dots,\ \rho_n\}$. The class
group of the latter variety is the cokernel of the function $M \to
\oplus_{i=1}^{n} {\Bbb Z} \cdot \rho_i$. The rows of the matrix for this
function are just
the transposed vectors $\eta_1, \dots,\ \eta_n$. So the class group is
determined using the methods from above.

As in
\cite[sequence (11)]{DFM:Cbg}, let
${\cal SF}$ denote the sheaf of support functions on $\Delta$, ${\cal M}$
the constant sheaf $M$ and ${\cal U}$ the sheaf kernel defined by the
sequence

\[
0 \to {\cal U} \to {\cal M} \to {\cal SF} \to 0
\]

\noindent
It follows from \cite{O:Cba} that $\Pic(X) \cong {\cal SF}(\Delta)/\im(M)$
and from
\cite[Theorem 1.a]{DFM:Cbg} that
\linebreak
 $H^2_{\acute{e}t}(K/X,G_m)$
 $\cong
\check{H}^1(\Delta,{\cal SF})$.
Let ${\cal L}$ denote the cokernel of the morphism ${\cal U} \to {\cal M}$ in
the category of presheaves. So ${\cal L}$ is a presheaf which is locally
isomorphic to ${\cal SF}$. That is, ${\cal SF}$ is the sheaf associated to
${\cal L}$. For any cone $\sigma \in \Delta$, let $\Delta(\sigma)$ denote
the subfan of $\Delta$ consisting of $\sigma$ and all of its faces. Then
${\cal L}(\Delta(\sigma)) = {\cal SF}(\Delta(\sigma))$ since support
funtions on a cone are linear.
Therefore we see that $\check{H}^p(\Delta,{\cal SF}) =
\check{H}^p(\Delta,{\cal L})$ for all $p$.
If $\{\sigma_1, \dots,\ \sigma_m\}$ are the maximal cones of $\Delta$,
and $\sigma_{ij}$ denotes $\sigma_i \cap \sigma_j$, then
the $\check C$ech complex
\[
0 \to \bigoplus_i {\cal L}(\Delta(\sigma_i))
\stackrel{\delta^0}{\longrightarrow}
 \bigoplus_{i<j} {\cal L}(\Delta(\sigma_{ij}))
\stackrel{\delta^1}{\longrightarrow}
 \bigoplus_{i<j<k} {\cal L}(\Delta(\sigma_{ijk}))
\]
can be used to compute the groups $H^p(\Delta,{\cal SF}) =
H^p(\Delta,{\cal L})$.
The  sequence
\[
0 \to {\cal U}(\Delta(\sigma)) \to  M \to {\cal L}(\Delta(\sigma)) \to 0
\]
is exact, so we see that ${\cal L}(\Delta(\sigma))$ is just the dual of the
group $N \cap {\Bbb R} \cdot \sigma$.
Given any cone $\sigma \in \Delta$, we find
${\cal L}(\Delta(\sigma))$ as follows.
First set up a matrix $S$ whose columns consist of those
$\eta_i$ such that $\rho_i \in \sigma$. A basis $L(\sigma_i)$
for ${\cal L}(\Delta(\sigma))$ is computed by methods of the introduction by
finding a basis for the smallest direct
summand of ${\Bbb Z}^r$ containing the column space of $S$.
For this and other
computations involving $L(\sigma)$
it is not necessary to distinguish between $N$ and its dual $M$.
If $\tau$ is a face of $\sigma$, there is a projection
${\cal L}(\Delta(\sigma)) \to {\cal L}(\Delta(\tau))$.
The matrix for this
projection corresponds to writing the elements in the
basis $L(\tau)$ in terms of the basis $L(\sigma)$ hence can be carried out
using the above algorithms. The ingredients for writing the matrix for
$\delta^0$ and $\delta^1$ are now available. Each is just a suitably
interpreted direct sum of projections of the form
${\cal L}(\Delta(\sigma)) \to {\cal L}(\Delta(\tau))$.

The kernel of $\delta^0$ is the group of support functions
${\cal SF}(\Delta)$.
As mentioned above, the Picard group is computed as
${\cal SF}(\Delta)/\im(M)$. The kernel of the map $M \to {\cal SF}(\Delta)$
is ${\cal U}(\Delta)$ which is just $H^0_{\acute{e}t}(X,G_m)/k^*$ hence we
can compute the group of units on $X$.
Again as mentioned above,
the relative Brauer group
$H^2_{\acute{e}t}(K/X,G_m)$ is
the first homology group $\ker(\delta^1)/\im(\delta^0)$.


\end{document}